\documentclass{amsart}

\newcommand{\PP }{{\mathbb P}}
\newcommand{\QQ }{{\mathbb Q}}
\newcommand{\CC }{{\mathbb C}}
\newcommand{\ZZ }{{\mathbb Z}}

\newcommand{\GG }{{\mathbb G}}

\newcommand{\ratmap}{--\!\!\!\!>}
\newcommand{\cal}{\mathcal}

\begin{document} 

\newcommand{\nt}{\noindent}
\newtheorem{theorem}{Theorem}
\newtheorem{conjecture}{Conjecture}

\title{Gromov-Witten Invariants for Abelian and Nonabelian Quotients}
\author{Aaron Bertram, Ionu\c t Ciocan-Fontanine and Bumsig Kim}
\maketitle
\section{Introduction.} Let $X$ be a
smooth projective variety over
${\CC}$ with the (linearized) action of a complex reductive group $G$, and let
$T \subset G$ be a maximal torus. In this setting, there are
two geometric invariant theory (GIT) quotients, $X//T$ and $X//G$, with 
a rational map $\Phi: X//T \ratmap X//G$ between them. We will further assume
that  ``stable $=$ semistable'' in the GIT and that all isotropy of stable
points is trivial, so $X//T$ and $X//G$ are smooth projective varieties, and
$\Phi$ is a $G/T$ fibration.

\medskip

Ellingsrud and Str\o mme \cite{ES} and Martin \cite{Mar} studied the relation
between the intersection theory of such quotients. In particular, there
is a lift of cohomology classes 
$\gamma \in H^*(X//G,{\QQ})$  to invariant classes $\widetilde 
\gamma \in H^*(X//T,{\QQ})^W;$
(for the action of the Weyl
group $W$), and Martin's {\it integration formula} relates the
Poincar\'e pairings:
$$\int_{X//G} \gamma \wedge \gamma' = \frac 1{|W|} \int_{X//T} 
(\widetilde \gamma
\wedge  \Delta^+) \wedge (\widetilde {\gamma'} \wedge \Delta^-)$$
where $\Delta^+ = \prod c_1(L_\alpha)$, product over the positive 
roots $\alpha$ (with 
line bundle $L_\alpha$), and 
$\Delta^-$ is the corresponding product for the negative roots.

\medskip

Gromov-Witten theory generalizes the intersection 
theory of a smooth projective variety $Y$ by means of intersection numbers on 
moduli spaces of maps from curves to $Y$.   
In \cite{HV}, Hori and Vafa made a ``physics'' conjecture relating Gromov-Witten 
theories of $X//G= \GG(s,n)$, the Grassmannian
of $s$-planes in $\CC^n$, and $X//T
=\PP^{n-1}\times\dots\times \PP^{n-1}=(\PP^{n-1})^s$ 
 (a mathematical 
version of the conjecture for genus zero curves 
was proved in our earlier paper \cite{BCK}), and
they suggested that their conjecture should extend to general flag manifolds.  

\medskip

In this paper, we contend that the appropriate generalized
context for the Hori-Vafa conjecture is that of nonabelian/abelian 
quotients described above and we 
state precise mathematical conjectures for
the genus zero theory.
In the second part of the paper, we prove the 
``$J$-function'' conjecture for flag manifolds.

\medskip

Given a smooth projective variety $Y$, 
a class $d \in H_2(Y,\ZZ)$, $n > 0$, cohomology 
classes $\gamma_1,...,\gamma_n \in H^{2*}(Y,\QQ)$ and 
$a_1,...,a_n \ge 0$, then
$\langle \tau_{a_1}(\gamma_1), \cdots, \tau_{a_n}(\gamma_n)\rangle_d \in\QQ$ 
denotes the associated (genus zero) Gromov-Witten invariant (see \S 3) 

\medskip

When a vector bundle $E$ on $Y$ is given, one can also define Gromov-Witten 
invariants twisted by a
multiplicative characteristic class
of $E$, see \cite{CG}. Our main conjecture (in a somewhat
imprecise form, see $(4.2)$ for the precise statement) can be viewed
as a ``quantum'' version of Martin's integration formula: 

\begin{conjecture} The genus zero Gromov-Witten invariants of $X//G$ are expressible
in terms of genus zero Gromov-Witten invariants of 
$X//T$ twisted by the Euler class
of the bundle $E=\oplus_{\alpha \in\{{\rm roots\; of}\; G\}}L_{\alpha}$.  
\end{conjecture}

The $n= 1$ invariants suffice for many applications to enumerative
geometry. Here, a cohomology-valued formal generating function of 
$(t_0,{\bf t})\in H^0(Y,\CC)\oplus
H^2(Y,\CC)$ (and an additional parameter $\hbar$) is formed:
$$J^Y(t_0,{\bf t},\hbar)=e^{t_0+{\bf t}/\hbar}
\sum_d e^{\int_d{\bf t}}J^Y_d(\hbar),$$
 with $J^Y_d(\hbar)$ defined by  
$\int_Y J^Y_d(\hbar) \wedge \gamma = \sum_{a=0}^\infty \hbar^{-a-2} 
\langle \tau_a(\gamma)\rangle_d.$
 
 \medskip

A special case of Conjecture 1 says then that
the $J$-function of $X//G$ can be calculated from the $J$-function
of $X//T$. 
Precisely, the relation is as follows:
Set
$$I_d(\hbar):=
\sum_{\widetilde d\mapsto d} \left(\prod_\alpha 
\frac{\prod_{k = -\infty}^{\widetilde d \cdot  c_1(L_\alpha)} (c_1(L_\alpha) + 
k\hbar)}{\prod_{k = -\infty}^{0} (c_1(L_\alpha) + 
k\hbar)}\right)J^{X//T}_{\widetilde d}(\hbar),$$
summed over all curve classes $\widetilde d\in H_2(X//T)$ lifting 
$d\in H_2(X//G)$ (see (4.1)),
and $$I(t_0,{\bf t},\hbar)=e^{t_0+{\bf t}/\hbar}\sum_d e^{\int_d{\bf t}}I_d(\hbar).$$ 

\begin{conjecture} $J^{X//G}(t_0,{\bf t},\hbar)$ is obtained 
from $I$ by an explicit change of variables (``mirror transformation'').
If $X//G$ is Fano of index $\geq 2$,  
$$J^{X//G}(t_0,{\bf t},\hbar)=I(t_0,{\bf t},\hbar).$$
\end{conjecture}

\nt {\bf Remarks:} Conjecture 2 resembles the quantum Lefschetz
theorem for ``concavex'' bundles (sums of ample and anti-ample line
bundles). Indeed, the
modification to $J_{\widetilde {d}}$ has exactly the form it would have if
$\oplus L_{\alpha}$ were concavex (which it isn't!) and $X//G$ were a complete 
intersection
in $X//T$ defined by the convex part of the bundle.  
See \cite{Kim2}, \cite{Kim3}
\cite{Lee}.
In fact, our most general conjectures (see (4.2) and (4.3)) 
include
a general version of quantum Lefschetz,  
involving an additional twist by Euler classes of
homogeneous vector bundles on $X//G$ and $X//T$.  

\medskip

Consider the flag manifold $F:=Fl(s_1,\dots s_l, n=s_{l+1})$
parametrizing flags:  $${\CC}^{s_1}\subset\dots\subset{\CC}^
{s_l}\subset {\CC}^n$$ 
and let $H_{i,j}, j = 1,...,s_i$ be Chern roots of the duals of the
universal bundles $S_i$:
$$S_1 \subset S_2 \subset \dots \subset S_l \subset S_{l+1} = 
{\CC}^n \otimes {\cal O}_{F}$$
and $d = (d_1,\dots, d_l)$ be the degree of a curve class, obtained by 
pairing 
with $c_1(S_i^\vee)$. 
Then Conjecture 2, together with Givental's formula for the $J$-function of 
the relevant toric variety gives the following closed formula for the 
$J$-function
of $F$ which we will prove and generalize to the other ``classical''
flag manifolds: 

\begin{theorem} For curve classes $d$ on the flag manifold $F$:
$$J^F_d(\hbar) = \sum_{\sum d_{i,j} = d_i} \prod_{i=1}^l 
\left( \prod_{1 \le j \ne j' \le s_i} \frac{\prod_{k = -\infty}^{d_{i,j} - d_{i,j'}}
(H_{i,j} - H_{i,j'} + k\hbar)}{\prod_{k = -\infty}^{0}
(H_{i,j} - H_{i,j'} + k\hbar)} \cdot\right.\ \ \ \ \ $$
$$\ \ \ \ \ \ \ \ \ \ \ \ \ \ \ \ \ \ 
\left.\prod_{1 \le j \le s_i, \ 1 \le j' \le s_{i+1}} 
\frac{\prod_{k = -\infty}^{0}
(H_{i,j} - H_{{i+1},j'} + k\hbar)}{\prod_{k = -\infty}^{d_{i,j} - d_{i+1,j'}}
(H_{i,j} - H_{{i+1},j'} + k\hbar)}   \right) 
$$
\end{theorem}
\medskip

\nt {\bf Remark.} With the same proof as in \cite{BCK} for 
Grassmannians, it follows from Theorem 1 and its generalization
 that Givental's $R$-Conjecture,
hence by \cite{Giv2} the Virasoro Conjecture, holds for classical generalized 
flag manifolds.
 
Earlier work on the $J$-function of type $A$ flag manifolds is contained in 
the recent paper \cite {LLY}, 
where a formula is given for $\int_Fe^{\bf t}J_d^F$,
but the problem of finding a closed formula for the function $J$ and 
a generalization of the
Hori-Vafa Conjecture is left open. The reader may want to compare our
Theorem 1 with the formula on page 39 of \cite{LLY}. Our proof of the
Theorem is a routine verification of our conjecture using Grothendieck
quot schemes.  The main point here is that once one has the correct
conjecture,  it is an easy matter to verify it. 
In particular, we do not use any of the results in \cite{LLY}.

\medskip

\nt {\bf Acknowledgments.}\quad Thanks are due to Alexander Givental, 
Dosang Joe, 
Yuan-Pin Lee, Jae-Suk Park, Ravi Vakil  and Youngho Woo for useful discussions. 
The final writing of this paper was done while the authors visited
National Center for Theoretical Sciences in Hsinchu, Taiwan, who we
thank for their hospitality and  wonderful working environment.  A. Bertram and I. Ciocan-Fontanine
were partially supported
by NSF grants DMS-0200895 and DMS-0303614, respectively. 
B. Kim was 
supported by KOSEF 1999-2-102-003-5, R02-2002-00-00134-0.

\begin{tableofcontents}

\end{tableofcontents}

\section {Classical intersection theory of GIT quotients.}

\nt {\bf (2.1)\quad Abelian and non-abelian quotients.} 

Let $X$ be a smooth projective variety with
fixed ample line bundle $\cal L$ and 
linearized action of $G$ as in the introduction.
 We denote the $G$-stable points by 
$X^s(G)$, respectively the $T$-stable points by $X^s(T)$, 
so that
$$X//G=X^s(G)/G\;\;\;\;  {\rm and}\;\;\;\; X//T=X^s(T)/T$$

The main example to have in mind is the situation considered in \cite{ES} where
there is a vector space $V$ on which $G$ acts linearly via a representation
$G\rightarrow GL(V)$ whose image contains the homotheties of $V$, so that
there is an induced action on $X=\PP (V)$. We adopt the notation in 
\cite{ES} in this case and write $V//G$ and $V//T$ for the quotients. 
Here it is clear that the unstable locus $X - X^s(G)$ has codimension $\ge 2$.
We will make this (mild) assumption as well in the general case. 

Under our hypotheses there is a diagram relating the two quotients:
\begin{equation}\begin{array}{ccccccc}
U& := & X^s(G)/T& \stackrel i\hookrightarrow & X^s(T)/T & =& X//T \\ \\
& & \downarrow \Phi  \\ \\
X//G & = & X^s(G)/G \\
\end{array}\label{diagram}\end{equation}

The map $i$ is an open immersion, while
$\Phi$ is a fibre bundle with fibre $G/T$, and can be further factored
as 
$$U\stackrel \varrho\longrightarrow Z\stackrel \eta\longrightarrow X//G$$
with $\eta$ a $G/B$-bundle and $\rho$ an (affine) $B/T$ bundle.
The diagram is constructed in detail in \cite{ES}, $\S 2$.

\medskip

Let $R$ be the root system of $G$ corresponding to the choice of the
torus $T$ and denote by $R^+$ the set of positive roots and by $R^-$
the set of negative roots. For each root $\alpha\in R$ there is an induced
line bundle $L_\alpha$ on $X//T$, coming from the canonical $1$-dimensional 
representation $\CC_\alpha$ of $T$ with weight $\alpha$. If $(\alpha,-\alpha)$
is a pair of opposite roots, then the corresponding line bundles are dual, 
and $c_1(L_\alpha)=-c_1(L_{-\alpha})$.

We denote
$$E^+:=\oplus_{\alpha\in R^+}L_{\alpha},\;\;\;\;
E^-:=\oplus_{\alpha\in R^-}L_{\alpha},\;\;\;\; E:=E^+\oplus E^- .$$
The Euler classes
$$\Delta =\Delta^+ ={\rm Euler}(E^+)=\prod_{\alpha\in R^+}c_1(L_{\alpha})$$ 
and 
$$\Delta^- =(-1)^{{\rm dim}(G/B)}\Delta
={\rm Euler}(E^-)=\prod_{\alpha\in R^-}c_1(L_{\alpha})$$
will play an important role in this paper.

Note that the Weyl group $W=N(T)/T$ of $G$ acts naturally on $X//T$, and 
therefore on the cohomology $H^*(X//T,\QQ)$. The classes $\Delta$ and 
$\Delta^-$ are $W$-anti-invariant, i.e. if $w$ is an element of $W$ of
length $\ell(w)$, then $w(\Delta)=(-1)^{\ell(w)}\Delta$.

\medskip

\nt {\bf (2.2)\quad Cohomology of $X//G$ versus cohomology of $X//T$.}

Unless mentioned otherwise, we will only 
consider cohomology with $\QQ$
coefficients. 
We recall some results of Ellingsrud-Str\o mme \cite{ES} 
and Martin \cite{Mar} relating the 
cohomology rings
$H^*(X//G)$ and $H^*(X//T)$. 

There are surjective Kirwan maps 
$$\kappa_G:H^*_G(X)\longrightarrow H^*(X//G)\;\;\;{\rm and}\;\;\;
\kappa_T:H^*_T(X)\longrightarrow H^*(X//T)$$
from equivariant cohomology of $X$ to the cohomology of the quotients, as
well as a natural restriction map
$$\tau_T^G:H^*_G(X)\longrightarrow H^*_T(X).$$
For cohomology classes $\gamma$ on $X//G$ and $\widetilde{\gamma}$ on $X//T$
we say that $\widetilde{\gamma}$ is a {\it lift} of $\gamma$ if they come 
from the same $G$--equivariant class on $X$:
$$\gamma=\kappa_G(\theta),\;\;\;\; 
\widetilde{\gamma}=\kappa_T(\tau_T^G(\theta))$$
for some $\theta\in H^*_G(X)$. From surjectivity of $\kappa_G$ it is clear 
that each $\gamma\in H^*(X//G)$ has a lift. 
Equivalently, one may define the notion of
lift using the maps in the basic diagram (\ref{diagram}) by the 
requirement that 
$$i^*(\widetilde{\gamma})=\Phi^*(\gamma).$$
This second description shows that $\widetilde{\gamma}$ may be taken
in $H^*(X//T)^W$. 
We repeat, for emphasis, Martin's {\it integration formula}:

\medskip

\nt {\bf (2.2.1)\quad Theorem (\cite {Mar}, Theorem B).} {\it With
notation as above,}
$$\int_{X//G}\gamma =\frac {1}{|W|}\int_{X//T}
\widetilde{\gamma} \wedge\Delta^+ \wedge
\Delta^-$$
{\it for any $\gamma \in H^*(X//G)$.}

\section {Genus zero Gromov-Witten theory.}

\nt {\bf (3.1)\quad Gromov-Witten invariants with descendents}

Let $Y$ be a smooth projective variety and let $d\in H_2(X,\ZZ)$ be fixed
curve class. The Kontsevich-Manin moduli stack
$\overline {M}_{0,n}(Y,d)$  of stable maps
of class $d$ from $n$-pointed nodal rational curves to $Y$ comes with
$n$ evaluation maps $ev_1,\dots ev_n$ to $Y$ (at the marked points).
The natural projection
$$\pi:\overline {M}_{0,n+1}(Y,d)\longrightarrow\overline {M}_{0,n}(Y,d)$$
given by forgetting the last marked point allows us to view
$\overline {M}_{0,n+1}(Y,d)$ as the universal curve. The map $\pi$ has $n$
sections $s_1,\dots ,s_n$ (corresponding to the marked points)
defining the {\it Witten cotangent
line bundles} $\cal L_i:=s_i^*(\omega_{\pi})$, with $\omega_{\pi}$ the relative
dualizing sheaf. It is customary to denote by $\psi_i$ the Chern class
$c_1(\cal L_i)$.

For cohomology classes $\gamma_1,...,\gamma_n \in H^{2*}(Y)$ and nonnegative
integers $a_1,...,a_n$, the associated genus zero Gromov-Witten invariant is
$$\langle \tau_{a_1}(\gamma_1), \cdots, \tau_{a_n}(\gamma_n)\rangle_d :=
\int_{[\overline {M}_{0,n}(Y,d)]^{\rm virt}}\wedge_{i=1}^n(\psi_i^{a_i}
\wedge ev_i^*(\gamma_i)),$$
where $[\overline {M}_{0,n}(Y,d)]^{\rm virt}$ is the virtual fundamental class
of \cite{LT}, \cite{BF}. This virtual class lives in
the Chow group
$A_{{\rm vdim}(\overline {M}_{0,n}(Y,d))}(\overline {M}_{0,n}(Y,d))$, with
$${\rm vdim}(\overline {M}_{0,n}(Y,d))=\int_d
c_1(\cal T_Y) + \mbox{dim}_{\CC}(Y) + n - 3$$ 
the virtual
dimension of $\overline {M}_{0,n}(Y,d)$, as given by the Riemann-Roch theorem.
The invariant is in general a rational number, and 
it vanishes unless 
$$\sum_{i=1}^n \frac{1}{2}{\rm deg}(\gamma_i) + a_i$$ equals the virtual
dimension. 

\medskip

\nt {\bf (3.2)\quad Gromov-Witten invariants twisted by the Euler class.} 
Assume now that on $Y$
we are given a vector bundle $E$. From:
$$\begin{array}{ccc} \overline {M}_{0,n+1}(Y,d)&\stackrel {e}
\longrightarrow & Y\\ \\ \downarrow \pi & &\\ \\
\overline {M}_{0,n}(Y,d) & & \\
\end{array}$$
with $e = ev_{n+1}$, we obtain
$$E_{n,d}=[R^0\pi_*e^*E]-[R^1\pi_*e^*E]$$
the push-forward of $e^*E$ in $K$-theory. Its virtual rank is
given by the Rieman-Roch formula: 
$${\rm vrk}(E_{n,d}):={\rm rk}(E)+ \int_dc_1(E).$$ 
Since $E_{n,d}$ can be obtained as the cohomology of a $2$-term complex
of vector bundles on $\overline {M}_{0,n}(Y,d)$ (see \cite{CG}), it has a 
well defined
``top Chern class'' $c_{\rm top}(E_{n,d}):=c_{\rm vrk}(E_{n,d})$. It is of
course zero if the virtual rank is negative, but this will not be the case
for any of the twisting bundles we employ.  

\medskip

We define the Gromov-Witten invariants of $Y$ twisted by the Euler class of $E$
by 
$$\langle \tau_{a_1}(\gamma_1), \cdots, \tau_{a_n}(\gamma_n)\rangle_{d,E} :=
\int_{[\overline {M}_{0,n}(Y,d)]^{\rm virt}}\wedge_{i=1}^n(\psi_i^{a_i}
\wedge ev_i^*(\gamma_i))\wedge c_{\rm top}(E_{n,d}).$$
These twisted invariants (in much greater generality) were studied by Coates
and Givental in \cite{CG}. Their definition of twisting by the
Euler class is formulated in a more general setting 
(essentially ``twisting by the Chern polynomial'') but can be seen to
specialize to the one given here for the cases we consider.

\section {Conjectural relations between $X//T$ and $X//G$}

Recall that on $X//T$ we have the bundle $E=\oplus_{\alpha\in R}L_{\alpha}$.
Martin's integration formula (2.2.1) may be viewed as expressing degree zero
Gromov-Witten invariants (i.e. usual intersection numbers) on $X//G$ in
terms of degree zero invariants on $X//T$ twisted by the Euler class of
the bundle $E$. We conjecture that this relation extends to stable maps
of higher degrees. 

\medskip

\nt {\bf (4.1)\quad Definition.} {\it For curve classes 
$d\in H_2(X//G, \ZZ)$
and $\widetilde{d}\in H_2(X//T,\ZZ)$ we say that $\widetilde{d}$ lifts $d$
(and write $\widetilde{d}\mapsto d$) if 
$$\int_dH=\int_{\widetilde{d}}\widetilde{H}$$
for every divisor class $H\in H^2(X//G,\QQ)$ with lift $\widetilde{H}\in
H^2(X//G,\QQ)^W$.}

\medskip

Since any two lifts agree when restricted to the
``$G$-stable locus'' $U\subset X//T$, and by assumption (see (2.1))
the complement of $U$ has codimension
at least 2, it follows that for divisor classes the 
$W$-invariant lifts
are unique, and the $\widetilde{d}$'s are indeed 
well-defined. (Note also that by taking $H$ to be ample in (4.1), with ample
$W$-invariant lift, we see that each $d$ has finitely many lifts.) 

As shown by the examples we treat in $\S 6$, it is actually useful to
introduce an additional twisting. Namely, consider a finite dimensional
linear representation $\cal V$ of $G$. It induces the homogeneous 
vector bundle $\cal V_G:=X^s(G)\times_G\cal V $ over $X//G$ and,
viewing $\cal V$ as a $T$-representation, the vector
bundle $\cal V_T:=X^s(T)\times_T\cal V $ over $X//T$. Since a
$T$-representation is completely reducible, $\cal V_T$ splits 
as a direct sum
of line bundles, which we will assume to be nef. 
Note that the Euler class of $\cal V_T$ is a lift of the
Euler class of $\cal V_G$.
With this definitions we can state:

\medskip

\nt {\bf (4.2)\quad Conjecture.} {\it 
Twisted genus zero Gromov-Witten invariants of $X//G$
and of $X//T$ are related by
$$\langle \tau_{a_1}(\gamma_1),\cdots,\tau_{a_n}(\gamma_n)\rangle_{d,\cal V_G}
=\frac{1}{|W|}
\sum_{\widetilde{d}\mapsto d}\langle \tau_{a_1}(\widetilde\gamma_1), 
\cdots, \tau_{a_n}
(\widetilde\gamma_n)\rangle_{\widetilde{d},E\oplus \cal V_T}, $$
where $\widetilde\gamma_i$ are lifts of $\gamma_i$.}

\medskip

In particular, when the extra twist by $\cal V$ is absent, we are
expressing the GW-invariants of $X//G$ in terms of invariants
of $X//T$ twisted by the Euler class of $E$, specializing to Martin's
formula in degree zero.
Next, we express this relationship in terms of Givental's $J$-functions.

By the work of Coates and Givental (\cite{CG}, Thm. 2 and Cor. 5),  
one can extract from (4.2) a conjectural
``Quantum Lefschetz formula'' calculating the $J$-function 
(on the ``big'' parameter space $H^{2*}(X//G)$!) of the nonabelian
quotient $X//G$
(or, more generally, the $J$-function of the theory on $X//G$ twisted by
$\cal V_G$)
in terms of the bundle $E$ (respectively, the
bundle $E\oplus \cal V_T$) and the $J$-function of 
the abelian quotient $X//T$. 
We state explicitly a weaker version involving the restriction of $J$
to the ``small'' parameter space $H^2(X//G)$. This restriction amounts
to considering (4.2) only for 1-point invariants.

We let $(t_0,{\bf t})$ denote a general element of $H^0(X//G,\CC)\oplus
H^2(X//G,\CC)$.
The $\cal V_G$-twisted $J$-function of $X//G$ is
$$J^{X//G,\cal V}(t_0,{\bf t},\hbar):=
e^{t_0+{\bf t}/\hbar}\sum_de^{\int_d{\bf t}}J_d^{X//G,\cal V_G}(\hbar),$$
with $J_d^{X//G,\cal V_G}(\hbar)$ defined by
$$\int_{X//G} J^{X//G,\cal V}_d(\hbar)\wedge\gamma\wedge c_{\rm top}(\cal V_G) 
=\sum_{a=0}^\infty \hbar^{-a-2} 
\langle \tau_a(\gamma)\rangle_{d,\cal V_G}.$$ 

Let $M_1,\dots ,M_r$ ($r$ is the dimension of the representation $\cal V$) 
denote the line bundle 
direct summands of the split bundle $\cal V_T$ on $X//T$ and assume all
the $M_i$ are nef, so that $\widetilde d\cdot c_1(M_i)\geq 0$ for every 
effective
curve class $\widetilde d$.
Define
$$I^{X//T,\cal V}(t_0,{\bf t},\hbar)=e^{t_0+{\bf t}/\hbar}\sum_de^{\int_d{\bf t}}
\sum_{\widetilde d\mapsto d} I_{\widetilde d}(\hbar),$$
where
$$I_{\widetilde d}(\hbar)=\prod_\alpha 
\frac{\prod_{k = -\infty}^{\widetilde d \cdot  c_1(L_\alpha)} (c_1(L_\alpha) + 
k\hbar)}{\prod_{k = -\infty}^{0} (c_1(L_\alpha) + 
k\hbar)}\prod_{i=1}^r 
\frac{\prod_{k = -\infty}^{\widetilde d \cdot  c_1(M_i)} (c_1(M_i) + 
k\hbar)}{\prod_{k = -\infty}^{0} (c_1(M_i) + 
k\hbar)}J^{X//T}_{\widetilde d}(\hbar).$$

In other words, $I^{X//T,\cal V}$ is 
obtained from the (untwisted) $J$ function of $X//T$ by
first introducing in each $J_{\widetilde d}$
the correcting classes determined by the bundles $E$ and $\cal V_T$, and
then specializing the parameter ${\bf t}$
to the subspace $H^2(X//T,\CC)^W$, which
we identify with $H^2(X//G,\CC)$. 

\medskip

\nt {\bf (4.2)\quad Conjecture.} {\it  There is an explicit change of variable
$(t_0,{\bf t})\mapsto f(t_0,{\bf t})$ 
such that
$$J^{X//G,\cal V}(t_0,{\bf t},\hbar)=I^{X//T,\cal V}(f(t_0,{\bf t}),\hbar).$$
}

\nt {\bf (4.3)\quad Remark.} As in the usual quantum Lefschetz hyperplane 
theorem (\cite {Giv1}, \cite{Kim2}, \cite{Lee}, \cite {CG}), the change of variable is 
read off
the asymptotics of the $1/\hbar$ expansion of the function $I$.
In particular,
if $X//T$ is Fano, and for every effective curve $C\subset X//T$ the 
intersection number $C\cdot (c_1(X//T)-\sum_{i=1}^rc_1(M_i))$ is at least 2,
then no change of variable is needed and we have the equality
$$J^{X//G,\cal V}(t_0,{\bf t},\hbar)=I^{X//T,\cal V}(t_0,{\bf t},\hbar).$$
This will be the case for all examples we treat in the rest of the paper.

Note that in general no analogue of the quantum Lefschetz
``correcting'' class is known when twisting by an indecomposable
vector bundle. 

\medskip

\section {Standard flag manifolds}

\nt {\bf (5.1)\quad The flag manifold as a GIT quotient.}
$F = Fl(s_1,\dots s_l, n=s_{l+1})$ from the introduction
is the GIT quotient: 
$$F = {\PP}(\oplus_{i=1}^l\mbox{Hom}({\CC}^{s_i},{\CC}^{s_{i+1}}))//G$$
by the action of $G=\prod_{i=1}^l\mbox{GL}(s_i,{\CC})$, where a matrix
$A\in \mbox{GL}(s_i,{\CC})$ acts on $\mbox{Hom}({\CC}^{s_i},{\CC}^{s_{i+1}})$
by left multiplication, and on $\mbox{Hom}({\CC}^{s_{i-1}},{\CC}^{s_i})$ by
right multiplication by $A^{-1}$. Stability here is the ordinary stability for
Grassmannians: an element
$x \in {\PP}(\oplus_{i=1}^l\mbox{Hom}({\CC}^{s_i},{\CC}^{s_{i+1}}))$ is stable if 
each of its ``coordinates'' $x_i \in \mbox{Hom}({\CC}^{s_i},{\CC}^{s_{i+1}}))$
is injective.
If $T \subset G$ is the product of the subgroups of diagonal
matrices, then the associated abelian quotient
$$Y:={\PP}(\oplus_{j=1}^l\mbox{Hom}({\CC}^{s_j},{\CC}^{s_{j+1}}))//T$$
is a toric variety. Corresponding to the description of the flag 
manifold as a tower of Grassmannian bundles:
$$\begin{array}{ccccc} \GG(s_i,s_{i+1}) & \rightarrow & Fl(s_i,s_{i+1},...,s_{l+1})
 & = & \GG(s_i,S_{i+1}) \\ & & \downarrow \\ & & Fl(s_{i+1},...,s_{l+1})
\end{array}$$
the toric variety $Y$ is a tower of product-of-projective-space bundles:
$$\begin{array}{ccccc} \prod^{s_i} {\PP}^{s_{i+1} - 1} & \rightarrow & 
Y_i & = & \PP(V_{i+1}) \times_{Y_{i+1}} \dots \times_{Y_{i+1}} \PP(V_{i+1})
\\ & & \downarrow \\ & & Y_{i+1}
\end{array}$$
with 
$$V_{i+1}=\bigoplus_{j=1}^{s_{i+1}}\cal O_{Y_{i+1}}(
\underbrace{0,\dots ,0,-1}_j,0,\dots ,0)$$
the  vector bundle
on $Y_{i+1}$ corresponding to $S_{i+1}$, with
$Y_{l+1} = \mbox{Spec}(\CC)$ and $Y_1 = Y$. 

\medskip

There is the additional right action of GL$(n,\CC)$ commuting with the action 
of $G$, descending to a (transitive) action on $F$, exhibiting $F$ as
a homogeneous space. When we speak of the equivariant cohomology
of $F$ and of $Y$, it will be with respect to the induced action of the 
maximal torus 
$T_n \subset GL(n,\CC)$.

\medskip 

The Chern roots $H_{i,j}, 1 \le i \le l, 1 \le j \le s_i$ of the 
introduction can 
now be 
viewed as Chern classes of the relative hyperplane classes of 
$Y_i \rightarrow Y_{i+1}$.
This way of representing the Chern roots is less efficient than the 
Chern classes on the full flag variety $Fl(1,2,...,n)$, since in that
case $H_{i,j} = H_{i+1,j}$ for each $j = 1,..., s_i$, whereas in this 
representation the $H_{i,j}$ are all distinct. The additional 
classes $H_{l+1,j}$ 
really appear extraneous, as they are Chern roots of the trivial bundle. But
they play an important role in the equivariant $J$-function of 
the toric variety, 
which we now describe, following Givental \cite{Giv1}.

\medskip

A toric variety $Y$ comes with a finite set $\{D_i\}$ of torus-invariant divisor
classes that generate $H^2(Y)$ additively and $H^*(Y)$ multiplicatively. 
In the description of $Y$ as a quotient ${\CC}^N//T$ for $T \subset (\CC ^*)^n$,
these are simply the first Chern classes $c_1(L_i)$, where $L_i$ are the $N$ 
line bundles determined by the coordinate lines of ${\CC}^N$. Moreover,
if $T \subset T' \subset ({\CC}^*)^n$, then the $T'/T$-equivariant 
cohomology ring of $Y$ is similarly generated by the equivariant divisor
classes determined
by the coordinate lines. 

\medskip

\nt {\bf (5.1.1)\quad Theorem (Givental, \cite{Giv1}).} {\it If $Y$ is a 
smooth toric variety
with the property that each curve class $d\in H_2(Y)$ satisfies 
$- \int_d K_Y 
\ge 2$, then the $J$-function of
$Y$ is given by the formula:}
$$J_d^Y(\hbar) = \prod_{i=1}^N \frac{\prod_{k = -\infty}^0 (D_i + k\hbar)}
{\prod_{k = -\infty}^{\int_d D_i}(D_i + k\hbar)}$$
{\it and the $T'/T$-equivariant $J$-function is given by the same formula,
with equivariant divisor classes.}

\medskip

\nt {\bf Note:} In case $Y = {\PP}^n$ or a product of projective spaces each 
$\int_d D_i \ge 0$, and the numerator can be factored out of the denominator, 
but in the 
general case this is the most convenient formulation of the $J$-function.

\medskip

For our toric variety, it is easy to see that the invariant divisors are:
$$H_{i,j} - H_{i+1,j'}, \ \ i = 1,...,l-1$$
in both the ordinary and equivariant case, and in addition
$H_{l,j} - H_{l+1,j'}$
which are equivariant if $H_{l+1,j'} = \lambda_{j'}$ in 
$H^*(BT_n) = {\QQ}[\lambda_1,...,\lambda_n]$
and ordinary if $H_{l+1,j'} = 0$.

\medskip

Thus, by Givental's theorem, we have:
$$J^{Y}_{\vec d}(\hbar) = \prod_{i=1}^l \ \ 
\prod_{1 \le j \le s_i, \ 1 \le j' \le s_{i+1}} \ \  
\frac{\prod_{k = -\infty}^{0}
(H_{i,j} - H_{{i+1},j'} + k\hbar)}{\prod_{k = -\infty}^{d_{i,j} - d_{i+1,j'}}
(H_{i,j} - H_{{i+1},j'} + k\hbar)}$$
where 
$$\vec d = (d_{1,1},...,d_{1,s_1},d_{2,1},...,d_{2,s_2},...,d_{l,1},...,d_{l,s_l})$$
is the general curve class (with, evidently, each $d_{l+1,j'} = 0$).

\medskip

It is also easy to see that the roots of $G$ give the divisor classes:
$$H_{i,j} - H_{i,j'}, \ \ i = 1,...,l,\ \ j\neq j'$$
so that our $J$-function conjecture becomes precisely Theorem 1. Again, we 
repeat for emphasis that the $T_n$-equivariant $J$-function is obtained by 
setting each $H_{l+1,j'} = \lambda_{j'}$ and the ``ordinary'' $J$-function is 
obtained by setting each $H_{l+1,j'} = 0$.

\medskip

\nt {\bf (5.2)\quad Conversion of the Formula.} 
In our previous paper \cite{BCK}, we 
proved that the $J$-function of the Grassmannian $\GG(s,n)$ is given by:
$$J^\GG(q,\hbar) = \sum q^d J^\GG_d(q,\hbar)$$
where
$$J^\GG_d(\hbar) = (-1)^{(s-1)d}\sum_{\stackrel{(d_1,\dots
,d_s)}{d_1+...+d_s=d}}
\frac{\prod_{1\leq j<j'\leq s}(H_j-H_{j'}+(d_j-d_{j'})\hbar)}
{\prod_{1\leq j<j'\leq s}(H_j-H_{j'})\prod_{j=1}^s\prod_{k=1}^{d_j}
(H_j+k\hbar)^n}$$
and $H_i$ are Chern roots of $S^\vee$, the dual of the universal subbundle.
We first rewrite the formula in a version more parallel
to the Givental formula for toric varieties:
$$
\sum_{\stackrel{(d_1,\dots ,d_s)}{d_1+...+d_s=d}}
\prod_{1\leq j\ne j'\leq s} \frac{\prod_{k =
-\infty}^{d_j-d_j'}(H_j-H_{j'}+k\hbar)} {\prod_{k =
-\infty}^{0}(H_j-H_{j'}+k\hbar)} \cdot \frac 1
{\prod_{j=1}^s\prod_{k=1}^{d_j} (H_j+k\hbar)^n}$$ 
and then it follows that the $J$-function for a 
{\it product} of Grassmannians:
$\prod \GG = \prod_{i=1}^l \GG(s_i,n)$
is given by:

$$(5.2.1)\;\;J^{\prod \GG}_{(d_1,...,d_l)}(\hbar) = 
\sum_{\vec d} \prod_{i=1}^l 
\left( \prod_{1 \le j \ne j' \le s_i} \frac{\prod_{k = -\infty}^{d_{i,j} - d_{i,j'}}
(H_{i,j} - H_{i,j'} + k\hbar)}{\prod_{k = -\infty}^{0}
(H_{i,j} - H_{i,j'} + k\hbar)} \times\right.  $$
$$\hskip 3in
\left.
\frac{1}{{\prod_{j=1}^{s_i}\prod_{k=1}^{d_{ij}}
(H_{i,j}+k\hbar)^n}}   \right)$$
where $\vec d$ is defined just as in the toric formula, with 
$\sum_{j=1}^{s_i} d_{i,j} = d_i$. Let
$S_i$ and $Q_i = \CC^n\otimes {\cal O}/S_i$ be the universal sub and quotient
$G$-bundles, thought of either on $\GG(s_i,n), \prod \GG(s_i,n)$ or $F$. Then
$F \subset \prod_{i=1}^l \GG(s_i,n)$
is a transverse zero section of the bundle
$\cal V = \bigoplus_{i=1}^{l-1} {\cal Hom}(S_i,Q_{i+1})$.
Of course, we may write:
$$0 \rightarrow \oplus_{i=1}^{l-1} {\cal Hom}(S_i,S_{i+1}) \rightarrow 
\oplus_{i=1}^{l-1}{\cal Hom}(S_i,{\bf C}^n) \rightarrow \cal V \rightarrow 0$$
and then the correction to the $d$th term of the $J$-function of 
$\prod \GG$ coming from 
twisting by $\cal V$, as
predicted by our general 
conjecture (4.2), is:
$$(5.2.2 )\;\;\;\;\;\;
\prod_{i=1}^{l-1}\left(\prod_{1 \le j \le s_i, \ 1 \le j' \le s_{i+1}} 
\frac{\prod_{k = -\infty}^{0}
(H_{i,j} - H_{{i+1},j'} + k\hbar)}{\prod_{k = -\infty}^{d_{i,j} - d_{i+1,j'}}
(H_{i,j} - H_{{i+1},j'} + k\hbar)}\times\right. $$
$$\hskip 3in
\left. {{\prod_{j=1}^{s_i}\prod_{k=1}^{d_{ij}}
(H_{i,j}+k\hbar)^n}}\right) $$
Or, in other words, our conjecture gives the same result whether we regard
the flag variety itself as a GIT quotient, or we think of it as 
the zero locus of a section of $\cal V$ in $\prod \GG$, regarded as a GIT quotient.
This is not a particularly deep check of the conjecture, but it is 
the latter point of view that we will use in the next section to prove the Theorem, 
as it generalizes immediately to the other classical Lie types.

\medskip

\section{Proof of Theorem 1 and its Generalizations}

\nt {\bf (6.1) A Simple J-function Lemma} 
The degree $d$ part of the ($T$-equivariant) $J$-function 
of $Y$ 
(with action of $T$) is given by the push-forward:
$$J^Y_d(\hbar) =  
e_*\left(\frac {[M^Y_d]}
{[M^Y_d/G^Y_d]}\right) =  
e_*\left(\frac {[M^Y_d]}
{\hbar(\hbar - \psi)}\right)$$
in ($T$-equivariant) cohomology, where we use the following conventions:

\medskip

$\bullet$ $M^Y_d = \overline M_{0,1}(Y,d)$ with virtual fundamental class
$[M^Y_d]$

\medskip

$\bullet$ $G^Y_d = \overline M_{0,0}(Y\times \PP^1,(d,1))$ with 
virtual fundamental class $[G^Y_d]$.
\medskip

$\bullet$ $\CC^*$ acts on $\PP^1$ by scaling $(x,y) \mapsto (tx,y)$ and 
$H^*(B\CC^*) = \QQ[\hbar]$.

\medskip

$\bullet$ Whenever $F\subset X$ is a fixed locus for an action of
$\CC^*$, then
$[F/X]$ denotes the Euler class of $F$ 
in $H^*_{\CC^*}(F) = H^*(F)[\hbar]$ (or $H^*_T(F)[\hbar]$), which is always invertible,
by the Atiyah-Bott localization theorem.

\medskip 

$\bullet$ $M^Y_d \subset G^Y_d$ is one of the fixed components for the induced 
$\CC^*$ action on $G^Y_d$. Specifically, it consists
of stable maps of curves with one component of class $(0,1)$ and 
the rest of the curve mapping to $0\in \PP^1$. 

It is a standard
fact in Gromov-Witten theory (see e.g. \cite{Ber}, \cite{Lee}) that:
$$[M^Y_d/G^Y_d] = \hbar(\hbar - \psi)$$ 

When $Y \subset \PP^n$, there is an equivariant
``map to the linear sigma model'' and diagram of
fixed components:
$$\begin{array}{rccccl} u:&  G^Y_d & \rightarrow & 
Y_d & \subset & 
\PP^n_d = \PP(\mbox{Hom}_d(\CC^2,\CC^{n+1})) \\
& \cup && \cup & & \cup \\
& M^Y_d &  \stackrel e \rightarrow &  Y & \subset & \PP^n
\end{array}$$
where $Y_d \subset \PP^n_d$ is defined by the equations 
induced from the equations of $Y$, and $\PP^n \subset \PP^n_d$
is the fixed locus arising from the ``zero'' line $\CC \subset \CC^2$ 
and inclusion $\PP^n = \PP(\mbox{Hom}_d(\CC,\CC^{n+1})) \subset \PP^n_d$.
Set-theoretically, $Y_d = \coprod _{e=0}^d \mbox{Map}_e(\PP^1,Y)\times
\PP^{d-e}$ so that if $Y$ is homogeneous (which will be our case),
then
$Y_d$ is the (singular!) closure of the (smooth) Hilbert scheme
$\mbox{Map}_d(\PP^1,Y) \subset \mbox{Map}_d(\PP^1,\PP^n) \subset \PP^n_d$
which is birational to the smooth compactification $G^Y_d$, and
$u_*[G^Y_d] = [Y_d]$.

\medskip

Now let 
$$i: X \stackrel k\hookrightarrow Y \stackrel j\hookrightarrow \PP^n$$
be $T$-equivariant embeddings for an action of $T$ on $\PP^n$ with 
isolated fixed points. Suppose $X$ and $Y$ are both  
homogeneous and
$\mbox{Map}_d(\PP^1,Y) \subset Q_d^Y$ 
is another smooth compactification with extended 
$\CC^*\times T$ action
and equivariant 
map $v: Q_d^Y\rightarrow \PP ^n_d$, and suppose there is an 
equivariant class
$[Q_d^X] \in 
A_*(Q^Y_d)$ such that:  
$$(\dag) \ \ v_*([Q_d^X])= [X_d] = u_*([G_d^X])$$  

Let $\alpha_F: F \hookrightarrow Q_d^Y$ be the union of fixed
components mapping to
$\PP^n$ by $v$. For ease of notation, we will pretend there
is only one component, writing, for example,
$\alpha_F^*[Q^X_d]/[F/Q^d_Y]$, when we really mean the sum 
$\sum \alpha_{F_k}^*[Q^X_d]/[F_k/Q^Y_d]$ over the 
components $F_k \subset F$. Let
$f: F \rightarrow \PP^n$ be the restriction of $v$. It follows
that $f$ factors through a map $g: F \rightarrow Y$, and  we get
the following diagram:
$$\begin{array}{cccccc} 
G_d^X & \stackrel u \rightarrow & \PP ^n_d & 
\stackrel v \leftarrow & Q_d^Y \\
& & \alpha_d\uparrow \ \ \ & &  \\ 
\alpha^X_d\uparrow \ \  & & \PP^n & & \alpha_F \uparrow \ \ \\
& & j \uparrow \ \ & \stackrel f\nwarrow \\
M_d^X & \rightarrow & Y & \stackrel g \leftarrow & F \\
& \stackrel e\searrow & k \uparrow \ \ \\
& & X
\end{array}$$

\nt {\bf (6.1.1) Lemma.} $[Q_d^X]$ computes $J^X_d$. First,
in $\CC^*$-equivariant cohomology:
$$\begin{array}{cccccccccccc}
i_*J^X_d & = & 
\frac{\alpha_d^*u_*([G_d^X])}{[\PP ^n/\PP ^n_d]}
& = & 
\frac{\alpha_d^*v_*([Q_d^X])}{[\PP ^n/\PP ^n_d]}
& = & f_* \frac{\alpha_F^*[Q_d^X]}{[ F/Q_d^Y]} 
& = & j_* g_* \frac{\alpha_F^* [Q_d^X]}{[ F/Q_d^Y]}
\end{array}$$
Next, in
$\CC^*\times T$-equivariant cohomology: 
$$\begin{array}{ccccccccccccc}
J^X_d &=& \frac{i^*i_*J_X^d}{[X/\PP ^n]} 
&=& \frac{1}{[X/\PP ^n]} i^*j_* g_* \frac{
\alpha_F^*[Q_d^X]}{[ F/Q_d^Y]}   
&=& \frac{1}{[X/Y]} k^*  g_*
(\frac{\alpha_F^*[Q_d^X]}{[ F/Q_d^Y]})
\end{array}$$

Finally, if $E$ is a $T$-equivariant vector bundle on $Y$, and $X$ is
the zero scheme of a section of $E$ transverse to the zero
section, then $[X/Y] = k^*c_{\rm top}(E)$, and:
$$J^X_d  =  k^* g_* (\frac{\alpha_F^*[Q_d^X]/g^*c_{\rm top}(E)}{[
F/Q_d^Y]})$$ and if (as will be our case) the right side is
well-defined as a $\CC^*$-equivariant cohomology class, then the
equality holds
in $\CC^*$-equivariant cohomology, taking limits. 

\medskip

\nt {\bf Remark:} The lemma is an easy exercise using the 
Atiyah-Bott localization theorem.
It immediately implies 
the quantum Lefschetz
theorem for quadrics, taking $Y = \PP^n$ and $Q_d^Y = \PP^n_d,$ with 
$[Q_d^X] = [X_d]$ defined by the $2d+1$ induced
quadratic equations.  Indeed, the obstruction to an easy proof of the
quantum Lefschetz theorem  for Fano complete intersections $X \subset
\PP^n$  is simply the fact that the linear sigma models $X_d \subset
\PP^n_d$ tend to have ``extra'' components aside from the closure of
$\mbox{Map}_d(\PP^1,X)$, making
property $(\dag)$ difficult to check. 
Such extra components do not exist when $X$ is a homogeneous  
space.

\medskip

\nt {\bf (6.2) Applying the Lemma.} Let 
$$Q^{\GG(s,n)}_d$$
be the Grothendieck quot scheme of vector bundle subsheaves
$K \subset \CC^n\otimes{\cal O}_{\PP^1}$
of degree $-d$ and rank $s$ on $\PP^1$. This is an alternative smooth compactification of 
the Hilbert space $\mbox{Map}_d(\PP^1,\GG(s,n))$ with universal 
sub and quotient sheaves:
$$0 \rightarrow {\cal K} \rightarrow 
\CC^n \otimes {\cal O}_{\PP^1\times Q^\GG_d} \rightarrow {\cal Q} \rightarrow 0$$
As we showed in our paper \cite{BCK}, under the natural map to the linear
sigma model:
$$\begin{array}{rccc} 
v: & Q^{\GG(s,n)}_d & \rightarrow & Q^{\GG(1,\binom{n}{s})}_d
= 
\PP^{\binom{n}{s} - 1}_d \\
& \cup & & \cup \\
& F & \stackrel g\rightarrow  &  \GG(s,n)
\end{array}$$
the components $F_{\vec d} \subset F$ are in bijection with the 
possible splitting types:
$$K = \oplus_{i=1}^s {\cal O}_{\PP^1}(-d_i)$$ 
Associated to each {\it splitting type} $\vec d = (d_1,d_2,...,d_s)$
(in non-decreasing order, with $d_1 = ... = d_{s_1}, d_{s_1+1} = ... = d_{s_2},...)$
is $F_{\vec d} \cong Fl(s_1,...,s_l=s,n)$ . The Euler classes to the $F_{\vec d}$ were
inverted and pushed forward in \cite{BCK}. This required 
a lemma of Brion (\cite{Bri}, see Lemma 1.4 of \cite{BCK}) and the following 
implication:

\medskip

\nt {\bf (6.2.1)} Let $P(x_1,...,x_s,d_1,...,d_s)$ be a polynomial
with the property that:
$$\prod_{1 \le i \le l} \prod_{\;\;\;\; 1 \le a < b \le s_i-s_{i-1}}
(x_{s_{i-1} + b} - x_{s_{i-1} + a}) \ \ \mbox{divides}\ \ 
P(x_1,...,x_s,\vec d)$$
for each splitting type $\vec d$ and such that each:
$$\frac {P(H_1,...,H_s,\vec d)}
{\prod_{i} \prod_{a < b}
(H_{s_{i-1} + b} - H_{s_{i-1} + a})}$$
represents a cohomology class in $F_{\vec d}$, where the $H_1,...,H_s$
are the Chern roots of $S$, arranged so that $H_{s_{i-1}+ 1},...,H_{s_i}$
are the Chern roots of $S_i/S_{i-1}$ on $F_{\vec d}$. Then:
$$\sum_{\vec d} g_*\left( \frac {P(H_1,...,H_s,\vec d)}
{\prod_{i} \prod_{a < b}
(H_{s_{i-1} + b} - H_{s_{i-1} + a})}\right) = 
\sum_{d_1 + ... + d_s = d}
\frac {P(H_1,...,H_s,d_1,...,d_s)}{\prod_{1 \le j < j' \le s} (H_{j'} -
H_j)}$$

 \nt {\bf Clarification:} The polynomial $P(x_1,...,x_r)$ in Lemma 1.4 
 of \cite{BCK} corresponds here to:
 $$\frac {P(x_1,...,x_s,\vec d)}
{\prod_{i} \prod_{a < b}
(x_{s_{i-1} + b} - x_{s_{i-1} + a})}$$
 
As we showed in \cite{BCK}:
$$\frac 1{[F_{\vec d}/Q_d^\GG]} = \frac{P_\hbar(H_1,...,H_s,\vec d)}
{\prod_{i} \prod_{a < b}
(H_{s_{i-1} + b} - H_{s_{i-1} + a})}$$
with
$$P_\hbar(x_1,...,x_s,d_1,...,d_s) = (-1)^{(s-1)(\sum d_i)}\frac{\prod_{1 \le j < j' \le s}
(x_{j'}- x_j + (d_{j'} - d_j)\hbar)}{\prod_{i=1}^s\prod_{k=1}^{d_i}
(x_i + k\hbar)^n}$$

Our point here is that (6.2.1) also applies to classes:
$$\frac{\alpha^*_{F_{\vec d}}[Q_d^X]/c_{\rm top}(E)}{[F_{\vec
d}/Q_d^\GG]}$$ (from(6.1.1)) describing the $J$-functions of
sub-homogeneous spaces. For example, consider the ``Lagrangian
Grassmannian:'' 
$$L\GG \subset \GG(n,2n) \subset \PP ^{\binom{2n}{n}-1}$$
obtained as the zero scheme of the section of $E = \wedge^2(S^\vee)$ over
$\GG(n,2n)$ induced from the standard algebraic symplectic form on
$\CC ^{2n}$. 

\nt {\bf Theorem 2.}
The $J$-function of the Lagrangian Grassmannian is given by:
\begin{eqnarray*}
J_d^{L\GG}  =  \sum _{
d_1 + ... + d_n =d} \left( \frac {\prod _{n\geq i > j\geq 1} \prod
_{k=0}^{d_i+d_j} (H_i +H_j+k\hbar)}{\prod
_{n\geq i > j\geq 1}(H_i+H_j)}\right) \cdot \\
\left(\frac{1}{\prod_{i=1}^n \prod _{k=1}^{d_i}(H_i
+ k\hbar )^{2n} }\cdot \prod _{n\geq i>j\geq 1}\frac{(H_i -H_j + (d_i-d_j)\hbar
)}{(H_i-H_j)}\right)
\end{eqnarray*}

{\bf Proof:} Consider the corresponding section of the vector bundle:
$$\pi_*\wedge^2 ({\cal K}^\vee)\ \ \mbox{on}\ \ 
Q_d^{\GG(n,2n)}$$
and its zero scheme, the ``Lagrangian Quot scheme'' $LQ_d^\GG \subset
Q_d^{\GG(n,2n)}$ (here $\pi$ 
is the projection $\PP^1\times Q_d^{\GG(n,2n)}\rightarrow\PP^1$). 
We apply (6.1.1) with:
$$[Q_d^{L\GG}] = c_{\rm top}(\pi_* (\wedge^2 {\cal K}^\vee))$$
noticing that $v_*[Q_d^{L\GG}] = u_*[G^{L\GG}_d]$ by virtue 
of the fact that $L\GG$ is homogeneous, and the zero section of 
$\pi_*\wedge^2{\cal K}^\vee$ is transverse along the 
Hilbert scheme of maps $\mbox{Map}_d(\PP^1,L\GG)$. So all that 
remains is the computation of $c_{\rm top}(\pi_*\wedge^2{\cal K})$ 
restricted to the fixed loci $F_{\vec d}$. But this is easily done. 
$F_{\vec d}$ consists of vector bundle subsheaves of splitting type $\vec d$, 
and there is a filtration of the restriction of ${\cal K}$ to 
$\PP^1\times F_{\vec d}$:
$${\cal K}_1 \subset {\cal K}_2 \subset ... \subset {\cal K}_l = {\cal
K} \subset {\cal O}_{\PP^1\times F_{\vec d}}$$  
with ${\cal K}_i/{\cal K}_{i-1} \cong \pi^*(S_i/S_{i-1})(-d_iD_0)$ where
$D_0 = 0\times F_{\vec d}$. From this, one computes
(as in \cite{BCK}, in the computations preceding Lemma 1.4):
$$c_{\rm top}(\pi_*\wedge^2{\cal K}^\vee)|_{F_{\vec d}} = 
Q_\hbar(H_1,...,H_n,\vec d)$$
where 
$$Q_\hbar(x_1,....,x_n,d_1,...,d_n) = \prod_{i>j}\prod_{k=0}^{d_i+d_j}
(x_i + x_j + k\hbar)$$
and then applying (6.2.1) to $P_\hbar \cdot
\left(Q_\hbar/\prod_{i>j}(x_i+x_j)\right)$ gives the Theorem, since 
$c_{\rm top}(E) = \prod_{i>j}(x_i+x_j)$. Note that evidently 
$\prod_{i>j}(x_i+x_j)$ divides $Q_\hbar$, hence the Theorem is valid
in $\CC^*$-equivariant cohomology.

\medskip

We next prove Theorem 1 in the same way, considering:
$$Fl(s_1,...,s_l,n) \subset \prod_{i=1}^l \GG(s_i,n) \subset \PP^N$$
with $E = \oplus_{i=1}^{l-1}{\cal Hom}(S_i,Q_{i+1})\cong \oplus_{i=1}^{l-1}
S_i^{\vee}\otimes Q_{i+1}$
(and $\prod \GG \subset \PP$ embedded by Pl\"ucker and Segre). Here the
zero scheme of the section of the bundle on the product of Quot schemes:
$$\pi_*\left(\oplus_{i=1}^{l-1} {\cal K}_i^\vee \otimes V/{\cal K}_{i+1}
\right)$$
is smooth and irreducible ($V$ is the trivial rank $n$ bundle). We will denote
the
 class of the zero section by $[Q^{H\GG}_{d_i}]$. It is the fundamental class
 of the
``hyperquot'' scheme of flags of vector bundle subsheaves of the trivial bundle on 
$\PP^1$:
$$K_1 \subset ... \subset K_l \subset \CC^n\otimes {\cal O}_{\PP^1}$$

Within the product of Quot schemes, the fixed components relevant to
our computations are the
products of flag manifolds $\prod_{i=1}^l F_{\vec d_i}$ indexed by multi-splitting types:
(Note: It is much harder to describe the fixed components of 
the hyperquot scheme. See \cite{LLY}). Now we proceed as in the proof of Theorem 2.
In this case, by (6.1.1) we need to compute the push-forward:
$$g_*\left(\sum_{\vec d_i} \alpha^*_{\prod F_{\vec d_i}}[Q^{H\GG}_{d_i}]/c_{\rm
top}(E)\right)$$  to the product of Grassmannians. This follows from an obvious
generalization of (6.2.1) to the product of Grassmannians. When applied to:
$$\prod P_\hbar(x_{i,1},...,x_{i,s_i},d_{i,1},...,d_{i,s_i})$$
we obtain the $J$-function for the product of Grassmannians. We need only, therefore,
to give the analogue of the polynomial $Q_\hbar$ from the proof of Theorem 2. For 
this, we notice that the sheaf ${\cal K}_i^\vee \otimes V/{\cal K}_{i+1}$ is not a vector bundle
(though its push-forward is a vector bundle). For this reason, it is more reasonable
to work with the  top Chern class via the exact sequence:
$$0 \rightarrow {\cal K}_i ^\vee \otimes {\cal K}_{i+1} 
\rightarrow {\cal K}_i^\vee \otimes V \rightarrow {\cal K}_i^\vee \otimes V/{\cal K}_{i+1} \rightarrow 0$$
and its pushforward to the product of quot schemes:
$$0 \rightarrow \pi_*({\cal K}_i ^\vee \otimes {\cal K}_{i+1} ) 
\rightarrow \pi_*( {\cal K}_i^\vee \otimes V) \rightarrow \pi_*({\cal K}_i^\vee \otimes V/{\cal K}_{i+1})
\rightarrow R^1\pi_*({\cal K}_i ^\vee \otimes {\cal K}_{i+1} ) \rightarrow 0$$
When the top Chern class is computed with this sequence, we obtain our desired 
polynomial
$Q_\hbar(x_{i,j},d_{i,j})$ which is divisible by $c_{\rm top}(E)$ and the quotient 
precisely evaluates to (5.2.2) when we set $x_{i,j} = H_{i,j}$. This proves Theorem 1.

\medskip

Finally, the reader may immediately generalize Theorems 1 and 2 to 
give the $J$-functions for all isotropic flag manifolds, i.e., the homogeneous
spaces $G/P$ for $G=SO_{2n+1}(\CC), Sp_{2n}(\CC), SO_{2n}(\CC)$ 
by realizing the
isotropic flag manifold inside the appropriate product of Grassmannians
as the zero locus 
of an appropriate homogeneous vector bundle.

\medskip

\noindent Department of Mathematics, University of Utah, Salt Lake
City, UT 84112,\\ {\it bertram@math.utah.edu}

\medskip

\noindent School of Mathematics, University of Minnesota,
Minneapolis MN, 55455,\\ {\it ciocan@math.umn.edu}

\medskip 

\nt School of Mathematics, Korea Institute for Advanced Study,
207-43 Cheongnyangni 2-dong, Dongdaemun-gu,
Seoul, 130-722, Korea,\\ {\it bumsig@kias.re.kr}

\end{document}